\documentclass[12pt]{article}

\usepackage{amssymb,latexsym,amsmath,amsfonts}

\oddsidemargin=\evensidemargin
\def\AGL{\mathop{\rm AGL}\nolimits}
\DeclareMathOperator{\AGaL}{A{\mathrm\Gamma}L}
\def\aut{\mathop{\rm Aut}\nolimits}

\def\cyc{\mathop{\rm Cyc}\nolimits}

\def\GL{\mathop{\rm GL}\nolimits}
\def\GR{\mathop{\rm GR}\nolimits}

\def\id{\mathop{\rm id}\nolimits}

\def\iso{\mathop{\rm Iso}\nolimits}

\def\rad{\mathop{\rm rad}\nolimits}
\def\rel{\mathop{\rm Rel}\nolimits}
\def\rk{\mathop{\rm rk}\nolimits}
\def\sym{\mathop{\rm Sym}\nolimits}

\def\A{{\cal A}}

\def\CC{{\cal C}}
\def\E{{\cal E}}

\def\H{{\cal H}}
\def\I{{\cal I}}

\def\M{{\cal M}}

\def\R{{\cal R}}
\def\X{{\cal S}}
\def\T{{\cal T}}
\def\U{{\cal U}}

\def\F{\mathbb F}

\def\ZZ{{\mathbb Z}}

\def\lg{\langle}

\def\ov{\overline}

\def\rg{\rangle}

\def\wh{\widehat}

\def\proof{{\bf Proof}.\ }
\def\bull{\vrule height .9ex width .8ex depth -.1ex }

\makeatletter
\renewcommand{\subsection}{\@startsection{subsection}{2}{0mm}{-2mm}{-2mm}
{\bf\normalsize}}

\def\sbsnt#1{\subsection{\hspace{-3mm}#1}}
\makeatother

\newtheorem{formula}{}[section]
\newtheorem{proposition}[formula]{Proposition}
\newtheorem{definition}{Definition}
\newtheorem{corollary}[formula]{Corollary}
\newtheorem{remark}[formula]{Remark}
\newtheorem{lemma}[formula]{Lemma}
\newtheorem{theorem}[formula]{Theorem}

\def\thrm{\begin{theorem}}
\def\thrml#1{\begin{theorem}\label{#1}}
\def\ethrm{\end{theorem}}
\def\prpstn{\begin{proposition}}
\def\prpstnl#1{\begin{proposition}\label{#1}}
\def\eprpstn{\end{proposition}}
\def\rmrk{\begin{remark}}
\def\rmrkl#1{\begin{remark}\label{#1}}
\def\ermrk{\end{remark}}
\def\dfntn{\begin{definition}}
\def\dfntnl#1{\begin{definition}\label{#1}}
\def\edfntn{\end{definition}}
\def\nmrt{\begin{enumerate}}
\def\enmrt{\end{enumerate}}
\def\tm#1{\item[{\rm (#1)}]}
\def\qtn{\begin{equation}}
\def\qtnl#1{\begin{equation}\label{#1}}
\def\eqtn{\end{equation}}
\def\lmm{\begin{lemma}}
\def\lmml#1{\begin{lemma}\label{#1}}
\def\elmm{\end{lemma}}
\def\crllr{\begin{corollary}}
\def\crllrl#1{\begin{corollary}\label{#1}}
\def\ecrllr{\end{corollary}}
\def\css{\begin{cases}}
\def\ecss{\end{cases}}

\begin{document}
\title{Normal cyclotomic schemes over a finite\\ commutative ring}
\author{
Sergei Evdokimov \\[-2pt]
\small Steklov Institute of Mathematics\\[-4pt]
\small at St. Petersburg \\[-4pt]
{\tt \small evdokim@pdmi.ras.ru }
\and
Ilia Ponomarenko\\[-2pt]
\small Steklov Institute of Mathematics\\[-4pt]
\small at St. Petersburg \\[-4pt]
{\tt \small inp@pdmi.ras.ru}
\thanks{The work was supported by RFFI Grants 05-01-00899
and NSH-4329.2006.1.}
}
\date{15.08.2006}
\maketitle

\begin{abstract}
We study cyclotomic association schemes over a finite commutative ring $R$ with identity.
The main interest for us is to identify the normal cyclotomic schemes $\CC$, i.e. those for which
$\aut(\CC)\le\AGaL_1(R)$. The problem is reduced to the case when the
ring $R$ is local in which a necessary condition of normality in terms of the subgroup of $R^\times$
defining $\CC$, is given. This condition  is proved to be
sufficient for a class of local rings including the Galois rings of odd characteristic.
\end{abstract}

\section{Introduction}\label{250706d}
Let $R$ be a finite commutative ring\footnote{Throughout the paper all rings are supposed to have identity.}
and $K$ a subgroup of its
multiplicative group $R^\times$. Denote by $\rel(K,R)$ the set of all binary relations on $R$ of the form
$\{(x,y)\in R\times R:\ y-x\in rK\}$, $r\in R$. Then the pair
\qtnl{100806b}
\cyc(K,R)=(R,\rel(K,R))
\eqtn
is an association scheme on $R$. We call it a {\it cyclotomic scheme over $R$} corresponding to the group~$K$.
Clearly, it is the scheme of 2-orbits of the group $\Gamma(K,R)=\{\gamma_{a,b}:\ a\in K, b\in R\}$ where $\gamma_{a,b}$ is
the permutation of the set $R$ taking $x$ to $ax+b$. In particular, it is a Cayley scheme over the additive
group $R^+$ of $R$ (see Subsection~\ref{220506a}) or a translation scheme in the sense of~\cite{BCN}.
Moreover, the multiplications by elements of $R^\times$ are Cayley isomorphisms of this scheme.

Cyclotomic schemes over a field were introduced by P.~Delsarte (1973) in connection with coding theory.
In \cite{EP01be} it was proved that any such scheme is uniquely determined up to isomorphism by its 3-dimensional
intersection numbers. Cyclotomic schemes over rings were introduced and studied in~\cite{GC92} within the framework
of the duality theory for association schemes. We also mention paper~\cite{IMY} where cyclotomic schemes over Galois
rings were used to construct amorphous association schemes. In the present paper we are interested in the
automorphism groups of cyclotomic schemes.

Historically, as the first result on the automorphism groups of cyclotomic schemes one should consider
the well-known Burnside theorem on permutation groups of prime degree. In fact, this theorem completely
determine the former groups for a prime field. In the case of an arbitrary finite field we have the following
result which is the interpretation of an old number-theoretical result from~\cite{M63} (see also~\cite[p.389]{BCN}).

\thrml{f080805b}
If $\CC$ is a cyclotomic scheme over a finite field $\F$, then $\aut(\CC)\le\AGaL_1(\F)$ whenever $\rk(\CC)>2$.\bull
\ethrm

For the cyclotomic schemes over the ring $\ZZ_n$ of integers modulo a positive integer $n$ the result of such a kind is not true.
Indeed, any such scheme being a Cayley scheme over a cyclic group $\ZZ_n^+$ can be treated up to language
as an S-ring over the same group. In accordance with \cite{LM,EP01be} each such S-ring can be constructed
from normal S-rings and S-rings of rank~2 by means of tensor products and generalized wreath products
(or wedge products in terms of \cite{LM}). Here normal S-rings are exactly those coming from
cyclotomic schemes $\CC$ such that $\aut(\CC)\le\AGaL_1(\ZZ_n)=\AGL_1(\ZZ_n)$. However, even among the S-rings
corresponding to cyclotomic schemes there exist non-normal ones.

The above discussion leads to the following definition which is central for this paper.

\dfntn
A cyclotomic scheme $\CC$ over a finite commutative ring $R$ is called {\it normal} if
$\aut(\CC)\le\AGaL_1(R)$.
\edfntn

The goal of the paper is to identify normal cyclotomic schemes.
Since any finite commutative ring is the direct product of local rings, the following theorem reduces the
general case to the local one (and, moreover, gives some product formula for two-points stabilizers of the
automorphism group).
Below for $R=\prod_iR_i$ we use the
following notation. For a cyclotomic scheme $\CC=\cyc(K,R)$ set $\CC_i=\cyc(K_i,R_i)$ where $K_i$ is defined from
the equality $\varphi_i(K_i)=K\cap \varphi_i(R_i^\times)$ with $\varphi_i$ being the monomorphism
of $R^\times_i$ to $R^\times$ such that the $j$th component of $\varphi_i(x)$ equals $x$ for $j=i$,
and equals $1_{R_j}$ for $j\ne i$.

\thrml{f090405a}
Let $R=\prod_iR_i$ be a finite commutative ring and $\CC$ a cyclotomic scheme over $R$. Then
\qtnl{f200106c}
\aut(\CC)_{u,v}=\prod_i\aut(\CC_i)_{u_i,v_i}
\eqtn
where $u=0_R$, $v=1_R$ and $u_i=0_{R_i}$, $v_i=1_{R_i}$ for all $i$.
In particular, $\CC$ is normal iff the scheme $\CC_i$ is normal for all $i$.
\ethrm

The following theorem gives a necessary condition for a cyclotomic scheme over a local ring to be normal. We
do not know any example when this condition is not sufficient. Below we set $I_0=\{x\in\rad(R):\, x\rad(R)=\{0\}\}$.

\thrml{f090405b}
Let a cyclotomic scheme $\cyc(K,R)$ over a finite local commutative ring~$R$ be normal. Suppose that
$K=K+I$ for some ideal $I$ of $R$. Then $I=\{0\}$ unless $q=2$ where $q$ is the order of the residue field of $R$.
Moreover, if $q=2$, then $I\subset I_0$.
\ethrm

Let $R$ be a local commutative ring. Given a group $K\le R^\times$ denote by $\I_K$ the set of all
ideals $I$ of $R$ such that $K+I=K$ or, equivalently, $1+I\subset K$. It is convenient for us to formulate the
following definition.

\dfntnl{100806c}
A group $K\le R^\times$ is called {\it pure} if the condition $I\in\I_K$ implies that $I=\{0\}$.
\edfntn

If $R$ is a field, then obviously any subgroup of
$R^\times$ is pure. Besides, Theorem~\ref{f090405b} implies that for $q>2$ the group $K$ is pure whenever the
scheme $\cyc(K,R)$ is normal. It turns out that for the Galois rings of odd characteristic other than fields
this necessary condition
of normality is also sufficient (as for the definition of a Galois ring see Section~\ref{f030406a}).

\thrml{f050106a}
Let $R$ be a Galois ring of odd characteristic other than a field. Then the scheme $\cyc(K,R)$ is normal
iff the group $K$ is pure.
\ethrm

Let $R=\GR(p^d,r)$ be a Galois ring of characteristic $p^d$ with the residue field of cardinality $q=p^r$
where $p$ is a prime. If $d>1$ and $p>2$ (the case of Theorem~\ref{f050106a}), then it is easy to see that a group
$K\le R^\times$ is pure iff it does not contain the group $1+p^{d-1}R$. On the other hand, if $d=1$, i.e.
$R=\F$ is a field of cardinality $q$, then the equality $\rk(\CC)=2$ implies that
$\aut(\CC)=\sym(\F)$. Besides, $\sym(\F)\le\AGaL_1(\F)$ iff $q\le 4$. Thus after combining
Theorems~\ref{f050106a} and~\ref{f080805b} we come to the following statement.

\thrml{250706a}
Let $R=\GR(p^d,r)$ with $p>2$. Then a cyclotomic scheme $\cyc(K,R)$ is normal
exactly in one of the following cases:
\nmrt
\tm{1} $d=1$ and either $(p,r)=(3,1)$ or $K\ne R^\times$,
\tm{2} $d>1$ and $K\not\ge 1+p^{d-1}R$.\bull
\enmrt
\ethrm

One of the ideas to prove the sufficiency in Theorem~\ref{f050106a} is to develop a reduction technique
for cyclotomic schemes over an arbitrary local ring $R$. For an ideal $I$ of $R$ the scheme $\cyc(\pi_I(K),R/I)$
where $\pi_I:R\to R/I$ is the natural epimorphism, can be treated as a factor-scheme of the scheme $\cyc(K,R)$
(see Subsection~\ref{230506a}). This simple observation is used in the proof of Theorem~\ref{f050805b} a
straightforward consequence of which is the following reduction statement. Below we set
$\pi_0=\pi_{I_0}$.

\thrml{f030406b}
Let $R$ be a finite local commutative ring, $\CC=\cyc(K,R)$ where $K\le R^\times$ is a pure group, and $\CC'=\cyc(K',R')$
where $K'=\pi_0(K)$ and $R'=R/I_0$. Then the scheme~$\CC$ is normal whenever so is the scheme~$\CC'$.\bull
\ethrm

Unfortunately, in the general case the group $K'$ is not pure (even if $R$ is a Galois ring of even
characteristic). So Theorem~\ref{f030406b} cannot be used for a direct inductive proof of the normality of
the scheme $\CC$. However, if $R$ is a Galois ring of odd characteristic, then this is true and Theorem~\ref{f050106a}
is reduced to the case $\rad(R)^2=\{0\}$. Thus, due to Theorem~\ref{f080805b} it suffices to prove the
following statement which is a special case of Theorem~\ref{f100406a}.

\thrml{f030406c}
Let $R$ be a finite local commutative ring other than a field for which $\rad(R)^2=\{0\}$.
Then the scheme $\cyc(K,R)$ is normal whenever the group $K$ is pure.\bull
\ethrm

Theorems~\ref{f050805b} and~\ref{f100406a} which are the origins of Theorems~\ref{f030406b} and \ref{f030406c} are
proved by using the S-ring technique. Namely, for a cyclotomic scheme $\CC$ over $R$ together with the ordinary (addition)
S-ring over $R^+$ corresponding to $\CC$ we consider its {\it multiplication} S-ring $\A$ over~$R^\times$ (see
Section~\ref{f120406c}). Everything is reduced to the case of a pure group $K\le\T\U_0$ where $\T$ is the
Teichm\"uller subgroup of $R^\times$ and $\U_0=1+I_0$. Then the group $\aut(\CC)_{u,v}$ acts faithfully on
$R^\times$ and the image of this action equals $\aut(\A)$. Moreover, in this case the S-ring $\A$ contains
the groups $\T$ and $\U=1+\rad(R)$, and becomes trivial after adding to it the cosets by any of these groups
(Section~\ref{010606a}). This enables us to prove that the group $\aut(\CC)$ normalizes the group $\AGL_1(R)$
(Theorems~\ref{f050406b} and~\ref{f180406a}). The latter means that $\aut(\CC)\le\AGaL_1(R)$ (Lemma~\ref{f250304a}),
i.e. the scheme $\CC$ is normal.
\vspace{2mm}

In fact, the developed technique permits us to obtain the following sufficient condition of normality for
an arbitrary finite local commutative ring $R$: the scheme $\cyc(K,R)$ is normal whenever the group $K$ is
strongly pure (Theorem~\ref{f050406a}). (Here we call a group $K\le R^\times$ {\it strongly pure} if it
is pure and the group $\pi_0(K)$ is strongly pure unless $R$ is a field.) It should be noted that this condition is not
necessary: one can prove that the cyclotomic scheme corresponding to the non-strongly pure group $K$
from the example in the beginning of Subsection~\ref{250706c} is normal.
\vspace{2mm}

In some cases one can say a little bit more on the automorphism group of a normal cyclotomic scheme $\CC=\cyc(K,R)$
where $R$ is a finite local commutative ring. For instance, if $K\le\T$ and $R$ is not a field, then
$$
\aut(\CC)\le\AGL_1(R)
$$
(statement (1) of Theorem~\ref{290506a}). This inclusion remains true also in some other
cases. In particular, this is so if the group $K$ is strongly pure, and either $K\le\U$ or
the residue filed~$\F$ of~$R$ is prime (statements (2) and (3) of Theorem~\ref{290506a}). The reason of this
is that in both cases the natural mapping
$$
\aut_\CC(R)\to\aut_\CC(\F)
$$
is a monomorphism and the group $\aut_\CC(\F)$ is trivial where by definition $\aut_\CC(R)$ (resp. $\aut_\CC(\F)$)
consists of all automorphisms of $R$ (resp. $\F$) that are automorphisms of $\CC$ (resp. the factor-scheme
of $\CC$ on $\F$) (see Theorem~\ref{f050406a}).
It should be noted that generally the kernel of the quotient homomorphism $\aut(R)\to\aut(\F)$ is not trivial.
For instance,  for $R=\F[X]/(X^2)$ the group
$\aut(R)$ is isomorphic to the semidirect product of~$R^\times$ by $\aut(\F)$ (indeed, the mapping
$a+b\pi\mapsto a^\sigma+b^\sigma\alpha\pi$
where $a, b\in{\F}$ and $\pi=X\mod X^2$, is an automorphism of~$R$ for any $\sigma\in\aut(\F)$ and $\alpha\in R^\times$).
\vspace{2mm}

All undefined terms and results concerning permutation groups can be found in~\cite{W64,W69,DM}. To make the paper
self-contained we cite the background on schemes and Schur rings in Section~\ref{010606b}. All necessary properties
of finite rings and cyclotomic schemes can be found in Section~\ref{f030406a}. The proofs of Theorems~\ref{f090405a}
and~\ref{f090405b} are contained in Section~\ref{010606c}, they are based on the ideas of~\cite{EP01be} where the case
$R=\ZZ_n$ was treated. The multiplication S-ring of a cyclotomic scheme is introduced and studied in
Sections~\ref{f120406c} and~\ref{010606a}. Section~\ref{f240406a} contains the proofs
of Theorems~\ref{f050805b}, \ref{f050406a}, and~\ref{f050106a}.
\vspace{2mm}

{\bf Notations.}
As usual by $\ZZ$ we denote the ring of integers.

For a ring~$R$ with identity we denote by $R^+$, $R^\times$ and $\rad(R)$ the additive and multiplicative groups
of~$R$ and the radical of~$R$ respectively.

Given groups $A\le R^\times$ and $B\le R^+$ with $AB=B$ we denote by $\Gamma(A,B)$ the group
$\{\gamma_{a,b}:\ a\in A,\ b\in B\}$ where $\gamma_{a,b}$ is the permutation of $R$ taking $x$ to $ax+b$. We omit
$B$ whenever $B=\{0\}$.

The group of all permutations of $V$ is denoted by $\sym(V)$.

Each permutation $f\in\sym(V)$ ($v\mapsto v^f$) naturally defines a permutation $R\mapsto R^f$ of the
set of all relations on~$V$. For an equivalence relation $E$ on a set $X\subset V$
such that $E^f=E$, the permutation $f$ induces a permutation $f^{X/E}\in\sym(X/E)$.
If all classes of $E$ are singletons, the set $X/E$ is identified with $X$.

For a group $G$ the permutation group on the set~$G$ defined by the right multiplications is
denoted by  $G_{right}$.

For $\Gamma\le\sym(V)$ and $X_1,\ldots,X_s\subset V$ we set $\Gamma_{X_1,\ldots,X_s}=\{\gamma\in\Gamma: X_i^\gamma=X_i$
for all~$i\}$. If $X_i=\{v_i\}$, the brackets are omitted. If the $X_i$'s are the classes of an
equivalence~$E$ on~$V$, we set $\Gamma_E=\Gamma_{X_1,\ldots,X_s}$.

\section{Finite commutative rings and cyclotomic schemes}\label{f030406a}

\sbsnt{Finite rings.}
It is well known (see e.g. \cite[Theorem~6.2]{MD74}) that any finite commutative ring
is the direct product of local rings. Let
$R$ be a finite local commutative ring. Then $R=\rad(R)\cup~R^\times$, the ideal $\rad(R)$ is maximal and the
characteristic of $R$ is a power of the characteristic of its residue field $\F=R/\rad(R)$. Moreover,
\qtnl{250706b}
R^\times=\T\times\U
\eqtn
where $\T$ is the Teichm\"uller group and $\U$ is the group of principal units. The groups $\T$ and $\U=1+\rad(R)$ are
a cyclic group of order $q-1$ and an abelian $p$-group respectively where $q$ and $p$ are the order and the characteristic
of~$\F$.

Let $I\subset\rad(R)$ be an ideal of $R$. Then the quotient ring $R/I$ is local and $(R/I)^\times=\pi_I(R^\times)$
where $\pi_I:R\to R/I$ is the natural epimorphism. Besides, the set $1+I$ is a subgroup
of $\U$. In particular, if $I\subset I_0$
then the mapping $r\mapsto 1+r$ induces an isomorphism of the additive group of $I$ onto $1+I$. Below we set $\U_0=1+I_0$.

The local ring $R$ is called {\it Galois} if $\rad(R)=pR$.\footnote{This is one of the equivalent definitions
given in \cite{MD74}.} Given positive integers $n,r$ there exists a unique up to isomorphism Galois ring of
characteristic $p^n$ with $q=p^r$; it is denoted by $\GR(p^n,r)$. We observe that $\GR(p,r)$ is a field of order
$p^r$ and $\GR(p^n,1)\cong\ZZ_{p^n}$. Each proper ideal of the  Galois ring $\GR(p^n,r)=R$ is of the form $p^iR$,
$i=1,\ldots,n$, and the quotient ring is isomorphic to $\GR(p^i,r)$. We also note that the homomorphism
$\aut(R)\to\aut(\F)$ induced by the epimorphism $\pi_{\rad(R)}$ is in fact an isomorphism (see~\cite{Wan}).

Generally, the structure of the group $\aut(R)$ (even in the local case) is unclear. Below we give a sufficient
condition for a permutation of $R$ to belong to this group.

\lmml{f250304a}
Let $R$ be a commutative ring and let a group $K\le R^\times$ generate $R^+$. Suppose that
$\gamma\in\sym(R)$ is such that
$$
0^\gamma=0,\quad 1^\gamma=1,\quad \gamma^{-1}\Gamma(K,R)\gamma=\Gamma(K,R).
$$
Then $\gamma\in\aut(R)$.
\elmm
\proof From the condition $\gamma^{-1}\Gamma(K,R)\gamma=\Gamma(K,R)$ it follows that given $(a,b)\in K\times R$ there exists
$(a_\gamma,b_\gamma)\in K\times R$ such that $\gamma^{-1}\gamma_{a,b}\gamma=\gamma_{a_\gamma,b_\gamma}$, or, equivalently,
\qtnl{f240204b}
(ax^{\gamma^{-1}}+b)^\gamma=a_\gamma x+b_\gamma,\quad x\in R.
\eqtn
Since $\gamma$ leaves fixed both $0$ and $1$, for $x=0$ this implies that $b_\gamma=b^\gamma$ for all $b\in R$,
whereas for
$(x,b)=(1,0)$ this implies that $a_\gamma=a^\gamma$ for all $a\in K$. Now, for $a=1$ and for
$b=0$ the equality (\ref{f240204b}) gives
\qtnl{f020805b}
(x+b)^\gamma=x^\gamma+b^\gamma,\quad (x,b)\in R\times R,
\quad\mbox{and}\quad
(ax)^\gamma=a^\gamma x^\gamma,\quad (a,x)\in K\times R
\eqtn
respectively. In particular, $\gamma\in\aut(R^+)$ and consequently (since $K$ generates $R^+$) the second equality holds for
all $a\in R$. Thus, $\gamma\in\aut(R)$.\bull
\vspace{2mm}

Lemma~\ref{f250304a} will be applied in Section~\ref{f240406a} to a local ring $R$ and $K=R^\times$. In this
case $\lg K\rg=R^+$ because any element of $\rad(R)=R\setminus R^\times$ is the difference of two units. One
more application is given by the following statement where $s=\gamma_{-1,1}$ is the involution
taking $x$ to~$-x+1$.

\crllrl{f180406c}
Let $\F$ be a field and $\gamma\in\sym(\F)$ a permutation leaving fixed both $0$ and~$1$.
Suppose that $\gamma$ normalizes both the groups $\Gamma(\F^\times)$ and $s\Gamma(\F^\times)s$.
Then $\gamma\in\aut(\F)$.
\ecrllr
\proof  A straightforward computation shows that
$\gamma_{a^{-1},0}s\gamma_{a,0}s=\gamma_{1,1-a}$ for all $a\in~\F^\times$. Then assuming (without loss of generality)
that $|\F|>2$ we see that the group $\lg\Gamma(\F^\times),s\Gamma(\F^\times)s\rg$ contains the group $\Gamma(1,\F^+)$
and hence equals $\Gamma(\F^\times,\F^+)$.
Since obviously $\lg\F^\times\rg=\F$, we are done by Lemma~\ref{f250304a} with $R=\F$ and $K=\F^\times$.\bull

\sbsnt{Cyclotomic schemes.}\label{230506a}
Let $\CC=\cyc(K,R)$ be a cyclotomic scheme over a finite commutative ring $R$ (see~(\ref{100806b})). Since
obviously each relation from $\rel(K,R)$ is $R^+_{right}$-invariant, $\CC$ is a Cayley
scheme over the group~$R^+$. The corresponding S-ring is called the {\it addition S-ring} of $\CC$. Each basic set of it is
of the form $rK$ where $r\in R$. It follows that any ideal $I$ of $R$ is an $\A$-subgroup (indeed,
$I=\bigcup_{r\in I}rK$). So due to the bijection between the sets $\H(\A)$ and $\E(\CC)$ (see Subsection~\ref{220506a})
we have the following statement.

\lmml{f2001006a}
For any ideal $I$ of $R$ the binary relation
$$
E(I)=\bigcup_{X\in R/I}X\times X
$$
belongs to the set $\E(\CC)$. In particular, the equivalence $E(I)$ is $\aut(\CC)$-invariant.\bull
\elmm
\vspace{2mm}

Since the set $\rel(K,R)$ is obviously $\AGL_1(R)$-invariant and the stabilizer of $u=0$ in $\AGL_1(R)$ equals
$\GL_1(R)$ we have
\qtnl{f130406b}
\AGL_1(R)\le\iso(\CC),\quad \GL_1(R)\le\iso(\CC_u)
\eqtn
where $\CC_u$ is the $u$-extension of $\CC$ (see Subsection~\ref{250706e}). The following easy statement gives a
simple criterion of normality.

\lmml{f020206a}
The scheme $\CC$ is normal iff $\aut(\CC)_{u,v}\le\aut(R)$ where $u=0$ and $v=1$.
\elmm
\proof The necessity follows from the obvious equality $\AGaL(R)_{u,v}=\aut(R)$. Conversely, by the orbit-stabilizer
property~\cite[Theorem~1.4A]{DM} we have
$$
[\aut(\CC): \aut(\CC)_{u,v}]=|R||K|=|\Gamma(K,R)|.
$$
Since $\Gamma(K,R)\le\aut(\CC)$, we conclude
that $\aut(\CC)=\aut(\CC)_{u,v}\Gamma(K,R)$ and the sufficiency follows.\bull
\vspace{2mm}

Now let the ring $R$ be local and $I\subset\rad(R)$ an ideal of $R$. Then $R/E(I)=R/I$, the equivalence
relation $E(I)$
is $\Gamma(K,R)$-invariant and $\Gamma(K,R)^{R/E(I)}=\Gamma(\pi_I(K),R/I)$. This implies that
\qtnl{250706f}
\cyc(K,R)_{R/E(I)}=\cyc(\pi_I(K),R/I),
\eqtn
i.e. the factor-scheme of $\CC$ modulo the equivalence $E(I)$ can naturally be treated as a
cyclotomic scheme over the ring $R/I$.

The following theorem on cyclotomic schemes with pure groups (see Definition~\ref{100806c}) will be used
in Section~\ref{f240406a}.

\thrml{f300805a}
Let $\CC=\cyc(K,R)$ be a cyclotomic scheme over a local commutative ring~$R$. If the group $K$ is pure, then
$$
\CC_{E_0}\ge\cyc(U_0,R)
$$
where $U_0=K\cap\U_0$ and $E_0=E(I_0)$.
\ethrm
\proof First we prove that if $S\in\rel(K,R)$ and $S_0\in\rel(U_0,R)$ are the relations
corresponding to sets $xK$ and $xU_0$ respectively, then
\qtnl{100806d}
S\cap((a+I_0)\times(b+I_0))=S_0\cap((a+I_0)\times(b+I_0)),\qquad a,b\in R,
\eqtn
whenever $x\in R^\times$ and the right-hand side is nonempty.
To do this let $(y,z)$ belong to the left-hand side. Then
$z-y\in (xK)\cap (b-a+I_0)$. On the other hand, due to the assumption there exists $(y_0,z_0)$ belonging
to the right-hand side. Then $z_0-y_0\in (xU_0)\cap (b-a+I_0)$. Thus, $(z-y)/(z_0-y_0)\in K\cap(1+I_0)=U_0$
and hence $z-y$ belongs to the right-hand side. The converse inclusion is obvious.

Denote by $\M$ the set of all relations from $\rel(U_0,R)$ corresponding to the sets $xU_0$ with $x\in R^\times$.
Then from (\ref{100806d}) it follows that $\M\subset\R^*(\CC_{E_0})$ and consequently $[\M]\le\CC_{E_0}$. Thus
it suffices to verify that $[\M]=\CC_{E_0}$, or, equivalently, that the addition S-ring $\A$ of the scheme $\cyc(U_0,R)$
is generated (as S-ring) by the sets $xU_0$, $x\in R^\times$. To do this we prove that
\qtnl{f290705c}
xU_0=\bigcap_{t\in\T}((x-t)U_0+tU_0),\quad x\in\rad(R).
\eqtn
Obviously, the left-hand side of (\ref{f290705c}) is contained in the right-hand side. Conversely, let
$t\in\T$ and $x\in\rad(R)$. Then since $U_0=1+H$ and $xH=\{0\}$ where $H$ is a subgroup of the group~$I_0$, we have
$$
(x-t)U_0+tU_0=(x-t)(1+H)+t(1+H)=x-t+tH+t+tH=x+tH.
$$
It follows that if $y$ belongs to the right-hand side of (\ref{f290705c}), then $y\in x+tH$ for all $t\in\T$.
On the other hand, $\bigcap_{t\in\T}tH=\{0\}$ by the purity of $U_0$. Thus, $y=x$ and we are done.\bull

\section{Proof of Theorems~\ref{f090405a} and~\ref{f090405b}}\label{010606c}

\sbsnt{Proof of Theorem~\ref{f090405a}.}
Set $\Gamma=\aut(\CC)$ and $\Gamma_i=\aut(\CC_i)$. To prove equality (\ref{f200106c}) first we
verify that $\prod_i(\Gamma_i)_{u_i,v_i}\le\Gamma_{u,v}$. To this end we observe that from the obvious inclusion
$\prod_iK_i\le K$ it follows that
$$
\prod_i\Gamma(K_i,R_i)\le\Gamma(K,R).
$$
So $\bigotimes_i\CC_i\ge\CC$ and hence $\prod_i\Gamma_i\le\Gamma$. Since obviously
$\prod_i(\Gamma_i)_{u_i,v_i}=(\prod_i\Gamma_i)_{u,v}$, we are done. To prove the converse inclusion
we observe that from
Lemma~\ref{f2001006a} with $I=R_i$ it follows that $R_i$ is a $\Gamma_u$-invariant set for
all $i$. For $\gamma\in\Gamma_u$ denote by $\gamma_i$ the restriction of $\gamma$ to $R_i$.
Then given an element $x=(\ldots,x_i,\ldots)$ of $R=\prod_iR_i$ we have
\qtnl{f200106d}
x^\gamma=(\ldots,x_i^{\gamma_i},\ldots).
\eqtn
Indeed, by Lemma~\ref{f2001006a} with $I=\prod_{j\ne i}R_j$ the equivalence $E(I)$ is $\Gamma_u$-invariant.
On the other hand, obviously each class of this equivalence contains a unique element of $R_i$.
Thus, the $i$th component of $x^\gamma$ equals $x_i^{\gamma_i}$ by the definition of $\gamma_i$.
Since $\Gamma_{u,v,R_i,1+R_i}=\Gamma_{u,v}$, from (\ref{f200106d}) it follows that
$$
\Gamma_{u,v}\le\prod_i(\Gamma_{u,v})^{R_i}=\prod_i(\Gamma_{u,v,R_i,1+R_i})^{R_i}\le
\prod_i((\Gamma_{R_i,1+R_i})^{R_i})_{u_i,v_i}.
$$
Thus the inclusion $\prod_i(\Gamma_i)_{u_i,v_i}\ge\Gamma_{u,v}$  and hence equality (\ref{f200106c})
are easy consequences of Lemma~\ref{f2001006b} below. Indeed, since the groups $\Gamma_i$ and $\Gamma(K_i,R_i)$
are 2-equivalent and the group $\Gamma_i$ is 2-closed, it follows that $(\Gamma_{R_i,1+R_i})^{R_i}\le\Gamma_i$.

\lmml{f2001006b}
For any $i$ the groups $(\Gamma_{R_i,{1+R_i}})^{R_i}$ and $\Gamma(K_i,R_i)$ are 2-equivalent.
\elmm
\proof Set $X=R_i$ and $Y=1+R_i$. Since $\Gamma(K_i,R)\le\Gamma$ and $\Gamma(K_i,R_i)=(\Gamma(K_i,R)_{X,Y})^X$,
it follows that $\Gamma(K_i,R_i)\le\Delta^X$ where $\Delta=\Gamma_{X,Y}$. So to prove the required statement it suffices to
verify that each 2-orbit of the group $\Delta^X$ is contained in some 2-orbit of the group $\Gamma(K_i,R_i)$,
or equivalently that each orbit of the group $(\Delta^X)_{u_i}=(\Delta_u)^X$ is contained in some orbit
of the group $\Gamma(K_i,R_i)_{u_i}=K_i$ (here we made used of the fact that the group $\Delta^X$ contains
a transitive subgroup $\Gamma(1,R_i)$). However, each orbit of the group $(\Delta_u)^X$ obviously meets
some orbit of the group $K_i$. So we only have to check that the latter orbit is $\Delta_u$-invariant. To do
this we need the following statement.

\lmml{f250106b}
For any $i$ and for any $a,r\in R$ with $r_j\in R^\times_j$ for all $j\ne i$, we have
$$
(a+rK)\cap R_i=a_i+r_is_iK_i
$$
for some $s\in K$, whenever the left-hand side set is nonempty.
\elmm
\proof From the definition of the monomorphism $\varphi_i$ it follows that $K=\bigcup_ssK'$
where $K'=\varphi_i(K_i)$ and $s$ runs over a full system of representatives of $K$ modulo $K'$.
Moreover, for all $s,t\in K$ we have
\qtnl{f030206a}
sK'=tK'\ \Leftrightarrow\ s_j=t_j\ \text{for all}\ j\ne i.
\eqtn
Besides,
\qtnl{f030206b}
a+rK=\bigcup_s(a+rsK')
\eqtn
for all $a,r\in R$. Suppose that the set $(a+rsK')\cap R_i$ is nonempty for some $a,r$ and $s$. Then $a_j+r_js_j=0$
for all $j\ne i$. So if $r$ is as in the lemma hypothesis, then the elements~$s_j$ for $j\ne i$, and hence
by (\ref{f030206a}) the coset~$sK'$, are uniquely determined by $a$ and $r$. Thus, in this case from
(\ref{f030206b}) it follows that
$$
(a+rK)\cap R_i=(a+rsK')\cap R_i=a_i+r_is_iK_i.
$$
Since the set $(a+rK)\cap R_i$ is nonempty iff the set $(a+rsK')\cap R_i$ is nonempty, we are done.\bull
\vspace{2mm}

Let us continue the proof of Lemma~\ref{f2001006b}.
First we observe that since $X+(v-v_i)=Y$, after translating by $v-v_i$ the equality of Lemma~\ref{f250106b}
with $a=v_i-v$ and $r=v-v_i$ we obtain that
$$
rK\cap Y=\{v-v_i\}.
$$
Denote by $S$ the basis relation of the scheme $\CC$ corresponding to $r$.
Since the sets $rK=S_{out}(u)$ and $Y$ are
$\Delta_u$-invariant, we conclude that so is the set $\{v-v_i\}$. Applying again Lemma~\ref{f250106b} for
$a=v-v_i$ and for $r$ such that the set $(v-v_i+rK)\cap X$ is nonempty and $r_j\in R_j^\times$ for all $j\ne i$,
we have
$$
(v-v_i+rK)\cap X=r_is_iK_i
$$
for some $s\in K$. Since the sets $v-v_i+rK=S_{out}(v-v_i)$ and $X$ are $\Delta_u$-invariant, we conclude that
so is the set $r_is_iK_i$. Besides, all the sets $r_is_iK_i$ cover $X$ when
$r$ runs over all elements of $R$ such that the set $(v-v_i+rK)\cap X$ is nonempty and $r_j\in R_j^\times$
for all $j\ne i$ (for instance, one can take $r_j=-1$ for $j\ne i$ and arbitrary $r_i\in R_i$).
Thus any orbit of the group $K_i$ is $\Delta_u$-invariant.\bull
\vspace{2mm}

Thus the first part of Theorem~\ref{f090405a} is proved. The second part follows from the first one and
Lemma~\ref{f020206a} applied to the schemes $\CC$ and $\CC_i$ for all~$i$.

\sbsnt{Proof of Theorem~\ref{f090405b}.}
Without loss of generality we assume that $R$ is not a field. Then the required statement is a straightforward consequence
of the following lemma the idea of the proof of which is taken from \cite[Subsection~5.3]{EP01be}.

\lmml{f080404a}
In the conditions of Theorem~\ref{f090405b} suppose that $R$ is not a field and $K+I=K$ for some nonzero ideal $I$ of
$R$. Then $|R/\rad(R)|=2$ and, moreover, $I\subset I_0$.
\elmm
\proof For each $k\in 1+I$ let us define a permutation $f_k$ of $R$ by
\qtnl{080806a}
x^{f_k}=
\css
kx, &\text{if $x\in\U$;}\\
 x, &\text{if $x\not\in\U$.}
\ecss
\eqtn
We show that $f_k\in\aut(\CC)$ where $\CC=\cyc(K,R)$. It suffices to verify that if $x-y\in rK$, then
$x^{f_k}-y^{f_k}\in rK$ for all $x,y,r\in R$. This is obvious for $x,y\not\in\U$, and follows from the
inclusion $1+I\le K$ for $x,y\in\U$. If $x\in\U$ and $y\not\in\U$, then
$$
x^{f_k}-y^{f_k}=kx-y=k(x-y)+(k-1)y\in rK+I=rK+rI=r(K+I)=rK
$$
and we are done. The other case is treated similarly.

From the normality of $\CC$ it follows that $x^{f_k}=ax^\sigma+b$ for some $a\in R^\times$, $b\in R$,
$\sigma\in\aut(R)$, and all $x\in R$. Since $0^{f_k}=0$ and $1^{f_k}=k$, we conclude that $b=0$ and $a=k$.
Thus, $x^{f_k}=kx^\sigma$, $x\in R$. Due to the choice of $k$ this implies that $\sigma$ leaves fixed
each element of $\U$ (and hence each element of $\rad(R)$) and each set $x+\rad(R)$. Since $\T^\sigma=\T$
and $|\T\cap(x+\rad(R))|=1$ for $x\in R^\times$, it follows that $\sigma$ leaves fixed each element of~$\T$.
Thus, $\sigma=\id_R$ and hence $x^{f_k}=kx$ for all $x\in R$. After comparing the latter equality with~(\ref{080806a})
we have
\qtnl{080806b}
(k-1)x=0,\qquad k\in 1+I,\ x\in R\setminus\U.
\eqtn
Since $\rad(R)\subset R\setminus\U$, we conclude that $k-1\in I_0$ for all $k$ and hence $I\subset I_0$. To complete
the proof suppose that $|R/\rad(R)|>2$. Then $R^\times\setminus\U\ne\emptyset$ and from (\ref{080806b}) it follows that $Ix=0$ for
some $x\in R^\times$. So $I=0$, which contradicts the choice of $I$.\bull

\section{Multiplication S-ring of a cyclotomic scheme}\label{f120406c}

Let $\CC=\cyc(K,R)$ be a cyclotomic scheme over a commutative ring $R$. Then from (\ref{f130406b}) it
follows that $\Gamma(R^\times)\le\iso(\CC_u)$ where $\CC_u$ is the $u$-extension of $\CC$ with $u=0$.
Since $\Delta(R^\times)$ is a relation of $\CC_u$ and the set
$R^\times$ is $\Gamma(R^\times)$-invariant, this implies that $R^\times_{right}=\Gamma(R^\times)^{R^\times}$
is a subgroup of $\iso((\CC_u)_{R^\times})$. So in accordance with Section~\ref{010606b} one can consider the scheme
$$
\CC'=((\CC_u)_{R^\times})^{R^\times_{right}}.
$$
Obviously, $R^\times_{right}\le\aut(\CC')$.
Thus, $\CC'$ is a Cayley scheme over the group $R^\times$. Denote by $\A=\A(K,R)$ the S-ring over $R^\times$
corresponding to the scheme $\CC'$.

\dfntn
The S-ring $\A$ is called the multiplication S-ring of the scheme~$\CC$.
\edfntn

The multiplication S-ring of a cyclotomic scheme over a field was introduced and studied in~\cite{EP01be}.

\thrml{f2130406c}
The set $\X^*(\A)$ contains both the sets $rK$ for all $r\in R^\times$ and the sets $(1+rK)\cap R^\times$
for all $r\in R$.
\ethrm
\proof Let $r\in R^\times$ and $C=rK$. Since each coset $C'\in R^\times/K$ is
the neighborhood of $u$ in the basis relation of $\CC$ corresponding to $C'$ , the set
$\Delta(C')$ is a relation of the scheme~$\CC_u$. Therefore, so is the relation $T$ defined by formula (\ref{f120506a})
with $G=R^\times$. Thus $C\in\X^*(\A)$ because the relation $T$ is $R^\times_{right}$-invariant
and $T_{out}(1)=C$.
To prove the second statement take $r\in R$ and set $X=(1+rK)\cap R^\times$. It is easily seen that
the smallest relation $S$ of $\CC_u$ containing $\{1\}\times X$ is a subset of $K\times R^\times$. Since all
relations of $\CC_u$ are $\Gamma(K)$-invariant, we see that $S_{out}(1)=X$. Besides, by the definition
of the scheme $\CC'$ the smallest relation~$S'$ of it containing $S$, is the union of all relations $Sr'=\{(sr',tr'):\ (s,t)\in S\}$
with $r'\in R^\times$. Since $S'_{out}(1)=S_{out}(1)=X$, it follows that $X\in\X^*(\A)$ and we are done.\bull
\vspace{2mm}

The following theorem establishes some connection between the automorphism groups of the scheme $\CC$
and the S-ring $\A$.

\thrml{f130406a}
In the above notation let $v$ be the identity of the ring $R$. Then
\nmrt
\tm{1} the mapping $f\mapsto f^{R^\times}$ induces a homomorphism from $\aut(\CC_{u,v})$ to $\aut(\A)$,
\tm{2} if $R$ is a field, then the mapping from statement (1) is an isomorphism.
\enmrt
\ethrm
\proof The set $R^\times$ being the neighborhood of $u$ in a relation of $\CC$, is $\aut(\CC_u)$-invariant. So the
mapping $f\mapsto f^{R^\times}$ induces a homomorphism from $\aut(\CC_u)$ to $\aut(\CC_u)^{R^\times}$. Besides,
$$
\aut(\CC_u)^{R^\times}\le \aut((\CC_u)_{R^\times})\le \aut(\CC').
$$
Since $\aut(\CC')_v=\aut(\A)$ by the definition of the latter group,
statement (1) follows. To prove statement (2) we note that the restriction homomorphism from $\aut(\CC_u)$ to
$\aut((\CC_u)_{R^\times})$ is an isomorphism inducing the isomorphism from $\aut(\CC_{u,v})$ to
$\aut((\CC_u)_{R^\times})_v$. On the other hand, by (\ref{f210406b}) with $\CC=(\CC_u)_{R^\times}$ and
$\Gamma=R^\times_{right}$ we have $\aut(\CC')=R^\times_{right}\aut((\CC_u)_{R^\times})$. So
$\aut(\A)=\aut(\CC')_v=\aut((\CC_u)_{R^\times})_v$ and we are done.\bull
\vspace{2mm}

In the general case the relationship between the groups $\aut(\CC_{u,v})$ and $\aut(\A)$ is unclear. However, we
have the following statement to be used in Section~\ref{f240406a}.

\thrml{f050406b}
Let $R$ be a local commutative ring, $\CC=\cyc(K,R)$ and $\A=\A(K,R)$. Then
the scheme $\CC_{u,v}$ is trivial whenever the S-ring $\A$ is trivial.
In particular, in this case $\aut(\CC)=\Gamma(K,R)$.
\ethrm
\proof Suppose that the S-ring $\A$ is trivial. This means that $\CC'$ is the scheme of 2-orbits of the
group $R^\times_{right}$ and consequently the scheme $(\CC')_v$ is trivial. So both the scheme
$((\CC_u)_{R^\times})_v$ and its extension $(\CC_{u,v})_{R^\times}$ are trivial too. On the other hand,
the permutation $s=\gamma_{-1,1}$ is an isomorphism of the scheme $\CC$ that interchanges $u$ and~$v$. So
$s\in\iso(\CC_{u,v})$. Thus, the scheme $(\CC_{u^s,v^s})_{(R^\times)^s}=(\CC_{v,u})_{1-R^\times}$ is trivial.
It follows that the restriction of $\CC_{u,v}$ to the set $R^\times\cup(1-R^\times)$ is trivial. Due to the locality
of the ring $R$ we have $R=R^\times\cup(1-R^\times)$. Thus the scheme $\CC_{u,v}$ is trivial.
The second part of the theorem follows from the first one and the proof
of Lemma~\ref{f020206a}.\bull

\section{Multiplication S-ring:\ pure case}\label{010606a}

In this section the multiplication S-ring of a cyclotomic scheme $\cyc(K,R)$ is studied
for a pure group $K\le\T\U_0$. First we rewrite the second half of the sets from Theorem~\ref{f2130406c}
in the multiplicative form.

\lmml{f230805c}
Let $R$ be a finite local commutative ring and $K=\T(1+H)\le R^\times$ with $H\le I_0$.
Then given $r=1+x\in\U$ we have
$$
(1+rK)\cap R^\times=\bigcup_{t\in\T,\ t\ne 1}t(1+z_{t,x}+\frac{t-1}{t}(H+x))
$$
where $z_{t,x}=y_tr$ with the element $y_t\in\rad(R)$ uniquely determined by the condition $1-t^{-1}+y_t\in\T$.
\elmm
\proof From (\ref{250706b}) it follows
\qtnl{f230805d}
(1+rK)\cap R^\times=\bigcup_{t\in\T,\ t\ne 1}((1+rK)\cap t\U).
\eqtn
Let $t\in\T$, $1+t\not\in\rad(R)$. Then $1+t=t'(1+y_{t'})$ for some $t'\in\T$. So
$$
(1+rK)\cap t\U=1+(1+x)t(1+H)=1+t(1+H+x)=(1+t)(1+\frac{t}{1+t}(H+x))=
$$
$$
t'(1+y_{t'})(1+\frac{t'(1+y_{t'})-1}{t'(1+y_{t'})}(H+x))=
t'(1+y_{t'}+\frac{t'(1+y_{t'})-1}{t'}(H+x))=
$$
$$
t'(1+y_{t'}+y_{t'}x+\frac{t'-1}{t'}(H+x))=
t'(1+z_{t',x}+\frac{t'-1}{t'}(H+x))
$$
(here $xH=y_{t'}H=\{0\}$ because $H\le I_0$).
To complete the proof, due to~(\ref{f230805d}) it suffices to note that $t'$ runs over the set $\T\setminus\{1\}$
when $t$ runs over the set $\T\setminus(-1+\rad(R))$.\bull

\thrml{f230805a}
Let $R$ be a finite local commutative ring, $K$ a subgroup of $R^\times$ and $\A$ an S-ring over $R^\times$ such that
$X(r)\in\X^*(\A)$ for all $r\in R$ where $X(r)=(1+rK)\cap R^\times$. Suppose that $K\le\T\U_0$ and the group $K$ is pure. Then
\nmrt
\tm{1} $\T,\U\in\H(\A)$,
\tm{2} the S-ring $\A$ is trivial whenever so is the S-ring~$\A_\T$ or the S-ring $\A_\U$.
\enmrt
\ethrm

Before proving Theorem~\ref{f230805a} we deduce from it an easy corollary to be used in the next section.

\thrml{f230805b}
Let $R$ be a finite local commutative ring and $\A=\A(K,R)$ where the group~$K$ is as in Theorem~\ref{f230805a}.
Then
\nmrt
\tm{1} $\T,\U\in\H(\A)$
\tm{2} if $H\in\{\T,\U\}$, then the S-ring generated by $\A$ and the cosets
of $R^\times$ by $H$ is trivial.
\enmrt
\ethrm
\proof Statement (1) immediately follows from statement (1) of Theorem~\ref{f230805a}.
because the S-ring $\A$ satisfies the hypothesis of this theorem (see Theorem~\ref{f2130406c}).
To prove statement~(2) denote by $\A'$ the S-ring generated by $\A$ and the cosets
of $R^\times$ by $H$. Since $\A'\ge\A$, the S-ring $\A'$ satisfies the hypothesis of Theorem~\ref{f230805a}.
So by statement (2) of that theorem it suffices to verify that the S-ring~$\A'_{H'}$ is trivial
where $H'=\U$ if $H=\T$, and $H'=\T$ if $H=\U$. However, this follows from the fact that $|H'\cap C|=1$
for any coset $C\subset R^\times$ by the group $H$.\bull
\vspace{2mm}

{\bf Proof of Theorem \ref{f230805a}.}
Since $\U$ is the complement of the set $\bigcup_{r\in R^\times} X(r)$ in $R^\times$, the second part of statement~(1)
follows. To prove the rest of the theorem we need the following lemma.

\lmml{f250406a}
In the condition of the theorem we have:
\nmrt
\tm{1} if $[tu_1]=[tu_2]$ for some generator $t$ of $\T$ where $u_1,u_2\in\U$, then $u_1=u_2$,
\tm{2} if $[t_1u]=[t_2u]$ for all $u\in\U$ where $t_1,t_2\in\T$, then $t_1=t_2$.
\enmrt
\elmm
\proof
Without loss of generality we assume that $\T\le K$. First, we prove statement~(1).
Any set $X\subset R^\times$ can uniquely be represented in the form $X=\bigcup_{t\in\T}tX_t$ where
$X_t\subset\U$ (see (\ref{250706b})). It follows that for any $\sigma\in\aut(\T)$ we have
\qtnl{f240805a}
X_t=(X^{\wh\sigma})_{t^\sigma},\quad t\in\T,
\eqtn
where $\wh\sigma$ is an automorphism of $R^\times$ such that $\wh\sigma^\T=\sigma$ and $\wh\sigma^\U=\id_\U$.
Since the group $\T$ is cyclic and its order is coprime to $|\U|$, the Chinese Remainder
Theorem implies that the automorphism $\wh\sigma$ is induced by raising to a power coprime to~$|R^\times|$.

Now let $X=[tu]$ where $t$ is a generator of $\T$ and $u\in\U$. Then obviously $u\in X_t$ and it suffices to verify that
\qtnl{f240805b}
X_t=\{u\}.
\eqtn
To do this we note that by the Schur theorem on multipliers $X^{\wh\sigma}=[t^\sigma u]$ for all $\sigma\in\aut(\T)$.
So $X^{\wh\sigma}\subset X(r_\sigma)$ for some $r_\sigma=1+x_\sigma$ with $x_\sigma\in\rad(R)$
(we made used of the fact that the latter set belongs to $\X^*(\A)$ and the union of all such sets equals
$R^\times\setminus\U$). Since $K\cap \U=1+H$ where $H\le I_0$, this implies by Lemma~\ref{f230805c} that
$$
(X^{\wh\sigma})_{t^\sigma}\subset 1+z_{t^\sigma,x_\sigma}+\frac{t^\sigma-1}{t^\sigma}(H+x_\sigma).
$$
However by (\ref{f240805a}) the element $u$ belongs to the left-hand side and so the right-hand side being
a coset by the group $\frac{t^\sigma-1}{t^\sigma}H$, equals
$u+\frac{t^\sigma-1}{t^\sigma}H$. Thus,
$$
X_t\subset\bigcap_{\sigma\in\aut(\T)}(u+\frac{t^\sigma-1}{t^\sigma}H)=u+H_0
$$
where $H_0$ is the intersection of all groups $\frac{t^\sigma-1}{t^\sigma}H$. To complete the proof of (\ref{f240805b})
we will show that $H_0=\{0\}$. Suppose on the contrary that there exists a nonzero $x\in H_0$. Then
$\frac{1}{1-t^\sigma}x\in H$ for all $\sigma\in\aut(\T)$. On the other hand, the following equality is true
$$
\frac{d}{1-t^d}=\sum_{i=0}^{d-1}\frac{1}{1-r^it},
$$
where $d$ is a proper divisor of $|\T|$ and $r$ is a generator of the subgroup of $\T$ of order $d$.
(It follows from the identity $\frac{f'(x)}{f(x)}=\sum_{i=0}^{d-1}\frac{1}{x-x_i}$ where
$f(x)=\prod_{i=0}^{d-1}(x-x_i)$ with $x_i\in R$, if we take $f(x)=x^d-1=\prod_{i=0}^{d-1}(x-r^i)$ and
put $x=t^{-1}$.) Since $d$ is coprime to the characteristic of $R$ and $r^it$ is a generator of $\T$ for all $i$,
the subgroup of $R^+$ generated by the set $\{\frac{1}{1-t^\sigma},\ \sigma\in\aut(\T)\}$ contains the set
$M=\{\frac{1}{1-t'}:\ t'\in\T\setminus\{1\}\}$. Thus, $Mx\subset H$. Since the image of $M$ in the residue
field $\F$ of $R$ equals $\F^\times\setminus\{1\}$ and $x\in I_0$, we have $Mx=Rx$. So
$1+I\subset 1+H\subset K$ where $I=Rx$, which contradicts the purity of $K$. This completes the proof
of statement~(1).

To prove statement (2) suppose that $[t_1u]=[t_2u]$ for all $u\in\U$ where $t_1,t_2\in\T$. By the second part of
statement (1) without loss of generality we can assume that $t_1\ne 1$ and $t_2\ne 1$. It suffices to verify that if
$t_1\ne t_2$, then there exists $r\in\U$ such that
\qtnl{f250805b}
t_1^{-1}(1+rK)\cap\U\ne t_2^{-1}(1+rK)\cap\U.
\eqtn
Indeed, in this case there exists an element $u$ belonging to the left-hand side but not belonging to the
right-hand side (or vice versa). Then $t_1u\in X(r)$ and $t_2u\not\in X(r)$.
Since $X(r)\in\X^*(\A)$, this implies that $[t_1u]\ne[t_2u]$, which contradicts the supposition.
To prove the required statement let $r=1+x\in\U$ be such that (\ref{f250805b}) becomes equality. Then
from Lemma~\ref{f230805c} it follows that
\qtnl{f250805c}
1+z_{t_1,x}+\frac{t_1-1}{t_1}(H+x)=1+z_{t_2,x}+\frac{t_2-1}{t_2}(H+x)
\eqtn
where $H$ is as above. Since the left-hand and right-hand sides are cosets by the groups
$\frac{t_1-1}{t_1}H$ and $\frac{t_2-1}{t_2}H$ respectively, these groups are equal.
Moreover, from (\ref{f250805c}) it follows that
$$
sx\in y+H'
$$
where $s=y_{t_1}-y_{t_2}+\frac{t_1-1}{t_1}-\frac{t_2-1}{t_2}$, $y=y_{t_2}-y_{t_1}$ (see Lemma~\ref{f230805c})
and $H'=\frac{t_1-1}{t_1}H=\frac{t_2-1}{t_2}H$.
We observe that $y\in\rad(R)$, and $s\in R^\times$ because $t_1\ne t_2$. It follows that
$x$ belongs to the coset
$$
C=s^{-1}y+s^{-1}H'\subset I_0.
$$
On the other hand, due to the purity of $K$ we have $H\ne I_0$
and hence $s^{-1}H'\ne I_0$. Thus, $C\varsubsetneqq I_0$ and inequality (\ref{f250805b}) is satisfied for any $r=1+x$ with
$x\in I_0\setminus C$.\bull

 To prove the first part of statement (1) take a generator
$t$ of the group $\T$. It suffices to verify that the set $X=[t]$ is contained in $\T$. To do this we observe
that from statement~(1) of Lemma~\ref{f250406a} it follows that $t^p\in X^{[p]}$ where $p$ is the characteristic
of the residue field of $R$ and $X^{[p]}$ is as in the Schur theorem on multipliers. So $t\in X'=(X^{[p]})^{\wh\sigma}$
where $\sigma$ is the automorphism of $\T$ inverse to raising to the $p$th degree and $\wh\sigma$ is the automorphism
of $R^\times$ defined above. Then by the Schur theorem on multipliers $X'\in\X^*(\A)$ and hence $X\subset X'$.
To complete the proof of statement (1) it suffices to observe that $|X'|\le|X|$ with the equality attained iff
$X\subset\T$.

To prove statement (2)
suppose that the S-ring $\A_\T$ is trivial. Let $u_1,u_2\in\U$, $u_1\ne u_2$. Then from statement (1) of
Lemma~\ref{f250406a} it follows that $[tu_1]\ne[tu_2]$ for some $t\in\T$. Since $[t]=\{t\}$, we have
$[tu_i]=[t][u_i]$ for $i=1,2$, whence it follows that $[u_1]\ne[u_2]$. Thus the S-ring $\A_\U$ is trivial and
consequently so is the S-ring $\A$. The second part of the statement is proved in a similar way by using
statement (2) of Lemma~\ref{f250406a}.\bull

\section{Strongly pure groups and proof of Theorem~\ref{f050106a}}\label{f240406a}

\sbsnt{Reduction statement.}
For a cyclotomic scheme $\CC$ over a ring $R$ we set
$$
\aut_\CC(R)=\aut(\CC)\cap\aut(R).
$$
If $I$ is an ideal
of $R$, we write $\aut_\CC(R/I)$ instead of $\aut_{\CC_{R/E(I)}}(R/I)$.

\thrml{f050805b}
Let $R$ be a finite local commutative ring and $\CC=\cyc(K,R)$ where $K\le R^\times$ is a
pure group. Then the natural mapping $\aut_\CC(R)\to\aut_{\CC}(R/I_0)$ is a
monomorphism. Moreover, $\aut(\CC)\le\AGaL_1(R)$ whenever $\aut(\CC)^{R/E_0}\le\AGaL_1(R/I_0)$
where $E_0=E_{I_0}$.
\ethrm
\proof We observe that the kernel of the homomorphism $f\mapsto f^{R/E_0}$ from $\aut(\CC)$ to
$\aut(\CC)^{R/E_0}$ coincides with the group $\aut(\CC)_{E_0}$. Moreover,
\qtnl{f170406a}
\aut(\CC)_{E_0}\le\AGL_1(R).
\eqtn
Indeed, $\CC_{E_0}\ge\cyc(U_0,R)$ by Theorem~\ref{f300805a}. So
$
\aut(\CC)_{E_0}=\aut(\CC_{E_0})\le\aut(\cyc(U_0,R)).
$
By Theorem~\ref{f050406b} it suffices to verify
that the S-ring $\A(U_0,R)$ is trivial. However, the latter immediately follows from Theorem~\ref{f230805b}.

Let $f\in\aut_\CC(R)$ be such that $f^{R/E_0}$ is identical. Then $f\in\aut(\CC_{E_0})$ and hence
$f\in\AGL_1(R)$ by (\ref{f170406a}). Since obviously $f$ leaves fixed the points $0$ and $1$, this
implies that $f=\id_R$. Thus the first statement of the theorem is proved.

To prove the second statement suppose that $\aut(\CC)^{R/E_0}\le\AGaL_1(R/I_0)$.
Then by Lemma~\ref{f250304a} with $K=R^\times$ it suffices to verify due to the locality of $R$
that the group $\Gamma=\Gamma(R^\times,R)$ is normalized by $\aut(\CC)$.
By~(\ref{f170406a}) all we need to prove is
$$
f^{-1}\Gamma f\subset\Gamma\aut(\CC_{E_0}),\quad f\in\aut(\CC).
$$
Take $\gamma\in\Gamma$ and $f\in\aut(\CC)$. Then
\qtnl{f050805c}
\ov{f^{-1}\gamma f}=\ov{f}^{-1}\ov{\gamma}\ov{f}\in \ov{f}^{-1}\ov{\Gamma}\ov{f}=\ov{\Gamma}
\eqtn
where the bar means factorization modulo $E_0$ (we made used of the fact that due to the assumption
$\ov{f}\in\ov{\aut(\CC)}\le\AGaL_1(R/I_0)$). On the other hand,
$f^{-1}\gamma f=\gamma(\gamma^{-1}f^{-1}\gamma)f=\gamma f_1$ where $f_1=(\gamma^{-1}f^{-1}\gamma)f$.
Since $\Gamma\le\iso(\CC)$, we have $f_1\in\aut(\CC)$. So from (\ref{f050805c}) it follows that
$$
\ov{f_1}\in\ov{\Gamma}\cap\ov{\aut(\CC)}=\Gamma(\ov{K},\ov R)
$$
where $\ov K=\pi_{0}(K)$ and $\ov R=R/I_0$. Due to the surjectivity of the natural homomorphism
$\Gamma(K,R)\to\Gamma(\ov{K},\ov{R})$, this
implies that there exists $\gamma_1\in\Gamma(K,R)$ such that $\ov{\gamma_1}=\ov{f_1}$. Thus,
$\gamma_1^{-1}f_1\in\aut(\CC_{E_0})$ and consequently
$$
f^{-1}\gamma f=(\gamma\gamma_1) (\gamma_1^{-1}f_1)\in \Gamma\aut(\CC_{E_0}).\bull
$$

\sbsnt{Strongly pure groups and normality.}\label{250706c}
We deduce Theorem~\ref{f050106a} from a more general result by using the notion of strong purity recursively
defined as follows. A group $K\le R^\times$ is called {\it strongly pure} if it is pure and the group
$\pi_0(K)$ is strongly pure unless $R$ is a field. Obviously, any strongly pure group is pure. The converse
statement is not true in the general case: a counterexample is given by $R=\F[X]/(X^n)$ where $\F$ is a finite field, $n\ge 3$, and
$K=1+\F x^{n-2}$ with $x=X\mod X^n$. However from the definition it immediately follows that it
is true whenever $\rad(R)^2=\{0\}$.

\thrml{f050406a}
Let $R$ be a finite local commutative ring other than a field. Then $\CC=\cyc(K,R)$ is a normal scheme whenever
$K$ is a strongly pure group. Moreover, in this case the natural mapping $\aut_\CC(R)\to\aut_{\CC}(\F)$ is a
monomorphism where $\F$ is the residue field of~$R$.
\ethrm
\proof With the help of Theorem~\ref{f050805b} applying inductively, the proof is reduced to the case $\rad(R)^2=\{0\}$.
Then the group~$K$ is pure. Thus the statement on monomorphism immediately follows from the first part of that theorem.
To complete the proof, assume (without loss of generality) that $K\ge\T$. Then  due to the second
part of the same theorem it suffices to verify that
\qtnl{f080805c}
\aut(\CC)^{\F}\le\AGaL_1(\F).
\eqtn
To do this we need the following lemma.

\lmml{250706g}
The groups $\aut(\CC_u)^{R^\times}$ and $\aut(\CC_v)^{1-R^\times}$ normalize the groups
$\Gamma(\T)^{R^\times}$ and $(s\Gamma(\T)s)^{1-R^\times}$ respectively where $u=0$, $v=1$ and $s=\gamma_{-1,1}$.
\elmm
\proof From statement (1) of Theorem~\ref{f230805b} it follows that $\T\in\H(\A)$
where $\A$ is the multiplication S-ring of the scheme $\CC$. So Theorem~\ref{f180406a} (with $H=\T$) implies that $\aut(\A)$ normalizes the group
$\lg \aut(\A'),\Gamma(\T)^{R^\times}\rg$. However, in our case the group $\aut(\A')$ is trivial
by statement (2) of Theorem~\ref{f230805b}. Therefore, $\aut(\A)$ normalizes $\Gamma(\T)^{R^\times}$. On the other hand, from
the definition of the S-ring $\A$ it follows that $\aut(\CC_u)^{R^\times}\le \Gamma(R^\times)^{R^\times}\aut(\A)$. Thus,
$\aut(\CC_u)^{R^\times}$ normalizes the group $\Gamma(\T)^{R^\times}$. To complete the proof we observe that obviously
$s$ is an isomorphism of $\CC$ that interchanges $u$ and $v$. Thus the group $\aut(\CC_v)^{1-R^\times}$ normalizes the group
$(s\Gamma(\T)s)^{1-R^\times}$.\bull
\vspace{2mm}

To prove (\ref{f080805c}) it suffices to show that if $\gamma\in\aut(\CC_{u,v})$, then $\gamma^\F\in\aut(\F)$.
However, from Lemma~\ref{250706g} it follows that the permutation
$\gamma^{X_0}=(\gamma^{R^\times})^{X_0}$ normalizes the group $\Gamma(\F^\times)^{X_0}$ whereas the permutation
$\gamma^{X_1}=(\gamma^{1-R^\times})^{X_1}$ normalizes the group $(s^\F\Gamma(\F^\times)s^\F)^{X_1}$
where $X_i=\F\setminus\{i\}$, $i\in\{0_\F,1_\F\}$.
Thus, $\gamma^\F$ normalizes both the groups $\Gamma(\F^\times)$ and $s^\F\Gamma(\F^\times)s^\F$ and we are done by
Corollary~\ref{f180406c}.\bull
\vspace{2mm}

Since a pure group $K\le\T\U_0$ is obviously strongly pure, from Theorem~\ref{f050406a} we obtain the
following statement.

\thrml{f100406a}
Let $R$ be a finite local commutative ring other than a field and $K$ a pure subgroup of $R^\times$.
Then the scheme $\cyc(K,R)$ is normal whenever $K\le\T\U_0$.\bull
\ethrm

We complete the subsection by giving a sufficient condition for the automorphism group of a cyclotomic
scheme to be a subgroup of $\AGL_1(R)$.

\thrml{290506a}
Let $\CC=\cyc(K,R)$ be a cyclotomic scheme over a finite local commutative ring $R$ other than a field. Then
$\aut(\CC)\le\AGL_1(R)$ whenever one of the following conditions is satisfied:
\nmrt
\tm{1} $K\le\T$,
\tm{2} $K$ is a strongly pure subgroup of $\U$,
\tm{3} the group $K$ is strongly pure and the residue field of $R$ is prime.
\enmrt
\ethrm
\proof To prove statement (1) we observe that from Theorem~\ref{f2130406c} it follows that the S-ring $\A(K,R)$
contains all elements of the set $R^\times/K$. So by statement (2) of Theorem~\ref{f230805b} this S-ring is trivial
and we are done by Theorem~\ref{f050406b}. To prove statements (2) and (3) we observe that from
the first part of Theorem~\ref{f050406a} it follows that $\aut(\CC)\le\AGaL_1(R)$. So by the second part
of this theorem it suffices to prove that the group $\aut_\CC(\F)$ is trivial where $\F$ is the residue
field of~$R$. However, for statement (2) this is clear because $\CC_\F=\cyc(\{1\},\F)$ (see (\ref{250706f}))
whereas for statement (3) because $\aut(\F)=\{1\}$.\bull

\sbsnt{Proof of Theorem \ref{f050106a}.}\label{010606d}
By Theorem~\ref{f050406a} the required statement is the consequence of the following theorem.

\thrml{f060206a}
For a Galois ring of odd characteristic any pure group is strongly pure.
\ethrm
\proof Let $R=\GR(p^n,r)$ where $p$ is odd. We can assume that $n>1$. Then the group~$\U$ is isomorphic
to the direct product of $r$ copies of a cyclic group of order $p^{n-1}$ (see \cite{MD74}). In particular, $\rk(\U)=r$ and
the maximal elementary abelian subgroup of $\U$ equals $\U_0$. Since the rank of
an abelian group does not increase by factorization, the theorem follows by induction from the lemma below.

\lmm
Under the above assumptions a group $K\le R^\times$ is pure iff $\rk(U)<r$ where $U=K\cap\U$.
\elmm
\proof Suppose that the group $K$ is not pure. Then $\U_0\le U$. On the other hand, since
$I_0=p^{n-1}R$, the group $\U_0$ is an elementary abelian of order $p^r$. Thus, $\rk(U)\ge\rk(\U_0)=r$.
Conversely, let $\rk(U)=r$. Then the maximal elementary abelian subgroup of $U$ is of order~$p^r$. So it
coincides with the maximal elementary abelian subgroup of $\U$. Since the latter subgroup equals $\U_0$, we are
done.\bull

\section{Appendix}\label{010606b}

\sbsnt{Schemes.}\label{250706e}
Let $V$ be a finite set and $\R$ a partition of $V^2$
closed with respect to transposition.\footnote{The elements of $\R$ are assumed to be nonempty.} Denote by $\R^*$ the set
of all unions of the elements of $\R$.
A pair $\CC=(V,\R)$ is called a {\it coherent configuration} or a {\it scheme} on $V$ if the diagonal $\Delta(V)$
of $V^2$ belongs to $\R^*$ and given  $R,S,T\in\R$, the number $|\{v\in V:\,(u,v)\in R,\ (v,w)\in S\}|$ does not depend on the
choice of $(u,w)\in T$. The elements of $V$, $\R=\R(\CC)$ and $\R^*=\R^*(\CC)$ are called the {\it points}, the {\it basis relations} and the {\it relations}
of~$\CC$, respectively. The number $\rk(\CC)=|\R|$ is called the {\it rank} of~$\CC$. We observe that
given $X,Y\subset V$ we have $X\times Y\in\R^*$ whenever $\Delta(X),\Delta(Y)\in\R^*$.

Two schemes are called {\it isomorphic} if there exists a bijection between their point sets preserving the basis
relations. Any such bijection is called an {\it isomorphism} of these schemes. The set of all isomorphisms
of a scheme $\CC$ is denoted by $\iso(\CC)$. This group contains a normal subgroup
$$
\aut(\CC)=\{f\in\sym(V):\ R^f=R,\ R\in\R\}
$$
called the {\it automorphism group} of~$\CC$. For a group $\Gamma\le\iso(\CC)$ we denote by $\CC^\Gamma$ the
scheme on the same point set the relations of which are exactly the elements of $\R^*$ invariant with respect
to $\Gamma$. In particular, if the scheme $\CC$ is {\it trivial}, i.e. $\R^*=2^{V^2}$, then $\iso(\CC)=\sym(V)$
and $\CC^\Gamma$ equals the {\it scheme of 2-orbits} of the group $\Gamma$. In the general case, it can be proved that
if $\Gamma$ acts regularly on the set $\{X\subset V:\ \Delta(X)\in\R\}$, then
\qtnl{f210406b}
\aut(\CC^\Gamma)=\Gamma\aut(\CC)
\eqtn
(see \cite[Theorem 2.2]{EP01be}).

Given a set $U\subset V$ denote by $\R_U$ the set of all nonempty relations $R_U=R\cap U^2$, $R\in\R$
(treated as relations on $U$). If $\Delta(U)\in\R^*$, then $\CC_U=(U,\R_U)$ is a scheme on $U$.
Clearly,
$$
\aut(\CC)^U\le\aut(\CC_U).
$$
Given an equivalence relation $E$ on $V$ denote by $\R_{V/E}$ the set of all relations
$$
R_{V/E}=\{(X,Y)\in (V/E)^2:\ R\cap (X\times Y)\ne\emptyset\}
$$
where $R\in\R$. If $E\in\R^*$,
then $\CC_{V/E}=(V/E,\R_{V/E})$ is a scheme on $V/E$. The set of all such $E$ is denoted by $\E=\E(\CC)$.
Clearly,
$$
\aut(\CC)^{V/E}\le\aut(\CC_{V/E}).
$$

The set of all schemes on $V$ is partially ordered by inclusion: namely, $\CC\le\CC'$ iff $\R^*\subset(\R')^*$.
The largest scheme is the trivial scheme on $V$ whereas the smallest one is the scheme of 2-orbits of the group $\sym(V)$.
For sets $\R_1,\ldots,\R_s$ of binary relations on~$V$ we denote by $[\R_1,\ldots,\R_s]$ the smallest scheme
$\CC$ on $V$ such that $\R_i\subset\R^*$ for all~$i$; we omit the braces if $\R_i=\{R_i\}$ and write $\CC_i$
instead of $\R_i$ if the latter is the set of basis relations of $\CC_i$. In particular, if $\CC$ is a scheme
on $V$ and $v_1,\ldots,v_s\in V$, then we set $\CC_{v_1,\ldots,v_s}=[\CC,\Delta(\{v_1\}),\ldots,\Delta(\{v_s\})]$.
One can see that
$$
\aut(\CC_{v_1,\ldots,v_s})=\aut(\CC)_{v_1,\ldots,v_s}
$$
where $\aut(\CC)_{v_1,\ldots,v_s}$ is the pointwise stabilizer of the set $\{v_1,\ldots,v_s\}$ in the group $\aut(\CC)$.
We will also use the following property of the {\it $v$-extension} $\CC_v$ of the scheme $\CC$ where $v\in V$: if $X=R_{out}(v)$
is the neighborhood of $v$ in a relation $R\in\R^*$, then $\Delta(X)$ is a relation of~$\CC_v$. Finally, if
$E$ is an equivalence relation on $V$, we set $\CC_E=[\CC,\{\Delta(X):\ X\in V/E\}]$. It immediately follows
that
\qtnl{f030805d}
\aut(\CC_E)\trianglelefteq\aut(\CC),\quad E\in\E(\CC).
\eqtn

\sbsnt{S-rings.}\label{220506a}
Let $G$ be a finite group. A subring~$\A$ of the group ring~$\ZZ[G]$ is called a {\it Schur ring} ({\it S-ring},
for short) over~$G$ if it has a (uniquely determined) $\ZZ$-base consisting of the elements
$\sum_{x\in X}x$ where $X$ runs over a family $\X=\X(\A)$ of pairwise disjoint nonempty subsets of~$G$ such that
$$
\{1\}\in\X,\quad
\bigcup_{X\in\X}X=G\quad
\textstyle{\rm and}\quad
X\in\X\ \Rightarrow\ X^{-1}\in\X.
$$
We call the elements of $\X$ {\it basic sets} of~$\A$ and denote by $\X^*=\X^*(\A)$ the set of all unions of them
and by~$\H=\H(\A)$ the set of all {\it $\A$-subgroups} of~$G$ (i.e. those belonging to~$\X^*$). The basic set of~$\A$
that contains $x\in G$ is denoted by $[x]$. The number $\rk(\A)=\dim_\ZZ(\A)$ is called the {\it rank} of~$\A$.
When $\rk(\A)=|G|$ (equivalently, $\A=\ZZ[G]$) we call the S-ring $\A$ {\it trivial}.
Given $H\in\H$ denote by $\A_H$ the S-ring over $H$ with $\X(\A_H)=\{X\in\X:\ X\subset H\}$.

The proof of the following theorem called the {\it Schur theorem on multipliers} can be found in~\cite{W64}. Below
given $X\subset G$ we set $X^{(m)}=\{x^m:\ x\in X\}$ for all $m\in\ZZ$, and
$X^{[p]}=\{x^p:\ x\in X,\ |xH\cap X|\not\equiv 0\pmod p\}$ for all prime $p$ where $H=\{g\in G:\ g^p=1\}$.

\thrm
Let $G$ be a finite abelian group and $\A$ an S-ring over $G$. Then for any $X\in\X(\A)$ the
following statements hold:
\nmrt
\tm{1} $X^{(m)}\in\X(\A)$ for any integer $m$ coprime to $|G|$,
\tm{2} $X^{[p]}\in\X^*(\A)$ for any prime $p$ dividing $|G|$.\bull
\enmrt
\ethrm

For a finite group $G$ denote by $\R(G)$ the set of all binary relations on $G$ that are invariant with respect
to $G_{right}$. Then the mapping
$$
2^G\to\R(G),\quad X\mapsto R_G(X)
$$
where $R_G(X)=\{(g,xg):\ g\in G,x\in X\}$, is a bijection. A straightforward computation shows that
if $H\triangleleft G$ and $C\in G/H$, then
\qtnl{f120506a}
R_G(C)=\bigcup_{C'\in G/H}C'\times CC'.\footnote{If $C=H$, then the normality of $H$ is not necessary.}
\eqtn
In particular, $R_G(H)$ is an equivalence relation on $G$.

Let $\A$ be an S-ring over the group $G$.
Then the pair $\CC=(G,\R)$ where $\R=R_G(\X)$, is a scheme on $G$ such that $G_{right}\le\aut(\CC)$.
Any scheme satisfying the latter condition is called a {\it Cayley scheme} on~$G$. In fact, the above
correspondence induces a bijection between the S-rings on $G$ and the Cayley schemes on~$G$ that preserves
the natural partial orders on these sets. Obviously,
$\R^*=R_G(\X^*)$ and $\E=R_G(\H)$. Moreover,
\qtnl{eoags}
\aut(\CC)=\aut(\A)\,G_{right}
\eqtn
where $\aut(\A)=\aut(\CC)_v$ with $v=1_G$.

\thrml{f180406a}
Let $\A$ be an S-ring over a group $G$ and $H$ a normal $\A$-subgroup of $G$. Then
$\aut(\A)$ normalizes the group $\lg\aut(\A'),H'\rg$ where $\A'$ is the S-ring over~$G$ generated
by $\A$ and the cosets of $G$ by $H$, and $H'$ is the subgroup of the group $G_{right}$ corresponding to
the multiplications by the elements of $H$.
\ethrm
\proof Let $\CC$ and $\CC'$ be the Cayley schemes over the group $G$ corresponding to the S-rings $\A$ and $\A'$
respectively. Then $\CC'=[\CC,\R]$ where $\R=\{R_G(C):\ C\in G/H\}$ (see (\ref{f120506a})). Let us show that
\qtnl{f190406a}
\CC_E=(\CC')_E
\eqtn
where $E=R_G(H)$. Indeed, since obviously $\Delta(C)$ is a relation of the scheme $\CC_E$, so is the
relation $R_G(C)$ for all $C\in G/H$. This implies that $\R\subset\R^*(\CC_E)$ and hence
$(\CC')_E\le\CC_E$. Since the converse inclusion is clear, equality (\ref{f190406a}) is proved. Next, we have
\qtnl{f190406b}
\aut((\CC')_E)=\aut(\A')H'.
\eqtn
Indeed, by definition $\aut(\CC')=\aut(\A')G_{right}$. Besides, due to the normality of $H$ we have
$(G_{right})_E=H'$.
Thus from the obvious equality $\aut(\A')_E=\aut(\A')$ it follows that
$$
\aut((\CC')_E)=\aut(\CC')_E=(\aut(\A')G_{right})_E=\aut(\A')(G_{right})_E=\aut(\A')H'
$$
whence (\ref{f190406b}) follows.

Since $H\in\H(\A)$, we have $E\in\E(\CC)$. So from (\ref{f030805d}) it follows that
$\aut(\CC_E)$ is a normal subgroup of $\aut(\CC)$. This implies that the group $\aut(\A)$ normalizes $\aut(\CC_E)$.
However $\aut(\CC_E)=\aut((\CC')_E)=\aut(\A')H'$ by (\ref{f190406a}) and (\ref{f190406b}), and we are done.\bull

\end{document}